\documentclass[12pt]{article}
\usepackage{tikz}
\usetikzlibrary{calc,through,backgrounds}
\usepackage{multicol}

\usepackage{bm,amssymb,amsthm,amsfonts,amsmath,latexsym,cite,color, psfrag,graphicx,ifpdf,pgfkeys,pgfopts,xcolor,extarrows}

\newtheorem{theo}{Theorem}[section]

\newtheorem{lem} [theo]{Lemma}

\newtheorem{prop}[theo]{Proposition}

\newtheorem{conj}[theo]{Conjecture}

\allowdisplaybreaks

\makeatletter \@addtoreset{equation}{section}
\@addtoreset{theo}{}\makeatother

\usepackage{hyperref}

\def\qed{\hfill \rule{4pt}{7pt}}
\def\pf{\noindent {\it Proof.} }

\begin{document}
\begin{center}
{\Large\bf Poincar\'e Polynomials of Odd Diagram Classes}

\vskip 6mm
{\small  Neil J.Y. Fan$^1$ and Peter L. Guo$^2$ }

\vskip 4mm
$^1$Department of Mathematics\\
Sichuan University, Chengdu, Sichuan 610064, P.R. China
\\[3mm]

$^{2}$Center for Combinatorics, LPMC\\
Nankai University,
Tianjin 300071,
P.R. China

\vskip 4mm

$^1$fan@scu.edu.cn, $^2$lguo@nankai.edu.cn
\end{center}

\noindent{A{\scriptsize BSTRACT}.}
An odd diagram class is a set of permutations with the same odd diagram.
 Brenti, Carnevale and Tenner showed that each odd diagram class
 is an interval in the Bruhat order. They  conjectured that
such intervals are rank-symmetric. In this paper, we present an algorithm
 to
 partition  an odd diagram class in a uniform manner.
 As an application, we obtain that the Poincar\'e polynomial of
an odd diagram class factors into polynomials of the form $1+t+\cdots+t^m$.
This in particular resolves  the conjecture of Brenti, Carnevale and Tenner.

\section{Introduction}

For two positive integers $a<b$, we use $[a,b]$ to represent the
interval $\{a, a+1, \ldots, b\}$.  Let $S_n$ denote the symmetric group of permutations of
$[n]:=[1,n]$. For   $w\in S_n$, we adopt the one-line notion, that is, we
write $w=w(1)w(2)\cdots w(n)$. An odd inversion of $w$ is an inversion   with the
parity condition, that is, a pair $(w(i), w(j))$ such that $1\leq i<j\leq n$,  $w(i)>w(j)$, and
$i\not\equiv j\ \pmod 2$. The odd length of $w$ is the  number of
  odd inversions of $w$. This statistic was introduced by Klopsch and   Voll \cite{KV} in
  their study of
  functions counting
non-degenerate flags in  formed spaces,   see also Brenti and  Carnevale \cite{BC}.

The odd diagram of $w$ is a diagram representation of its odd inversions,
which can be viewed as an odd analogue  of the classical Rothe diagram of $w$. Specifically,
the odd diagram $D_o(w)$ of $w$ is a subset of boxes in an $n\times n$ square grid
defined by
\[D_o(w)=\{(i,j)\colon w_i>j,\, i<w^{-1}(j),\, i\not\equiv w^{-1}(j) \pmod 2\},\]
where $w^{-1}$ is the inverse of $w$.
Here, we use the matrix coordinates, and use $(i,j)$ to denote the box in row $i$ and column $j$.
A subdiagram $D$ is called an odd diagram if there exists $w\in S_n$
such that $D_o(w)=D$. For an odd diagram $D$, let  $\mathrm{Perm}(D)$
denote the odd diagram class of $D$, namely,
\[\mathrm{Perm}(D)=\{w\in S_n\colon D_o(w)=D\}.\]

Brenti,   Carnevale and  Tenner \cite{BCT} proved that
 odd diagram classes partition the symmetric group in
  an extremely  pleasant way.

\begin{theo}[\mdseries{Brenti--Carnevale--Tenner \cite[Theorem B]{BCT}}]\label{interval}
Each odd diagram class    is an  interval in the Bruhat order.
\end{theo}

They  conjectured that   $\mathrm{Perm}(D)$
satisfies   a stronger symmetry   property.

\begin{conj}[\mdseries{Brenti--Carnevale--Tenner \cite[Conjecture 6.12]{BCT}}]\label{CJ}
Each odd diagram class   is rank-symmetric in the Bruhat order.
\end{conj}

For a Bruhat interval $[u, v]$ in $S_n$,
the associated Poincar\'e polynomial $P_{u, v}(t)$
is  the rank  generating function:
\[P_{u, v }(t)=t^{-\ell(u)}\sum_{u\leq w\leq v} t^{\ell(w)},\]
where $\ell(w)$ is the Coxeter length of $w$.
In the case when $[u, v]$ is a lower interval $[e, w]$, $P_w(t^2):=P_{e, w}(t^2)$
specifies to the Poincar\'e polynomial of the cohomology ring of the Schubert
variety $X_w$ indexed by $w$.

In this paper, we prove that the Poincar\'e polynomial of $\mathrm{Perm}(D)$
admits the following factorization.

\begin{theo}\label{main}
The Poincar\'e polynomial of an odd diagram class  can be expressed as a product
of factors of the form $1+t+\cdots+t^m$.
\end{theo}

A polynomial $f(t)=a_0+a_1t+\cdots+a_dt^d$ of degree $d$ is called
palindromic if
\[t^d\, f(t^{-1})=f(t).\]
Clearly, a Bruhat interval is rank-symmetric if and only if the associated
Poincar\'e polynomial
is  palindromic.
Theorem \ref{main} obviously implies that  the Poincar\'e polynomial of $\mathrm{Perm}(D)$
is  palindromic, thus confirming  Conjecture \ref{CJ}.

\noindent
{\bf Remark.} It is well known that the following for  $w\in S_n$ are equivalent:
\begin{itemize}

\item[(1)] the Schubert variety $X_w$   is smooth;

\item[(2)] the Poincar\'e polynomial $P_w(t)$   is  palindromic;

\item[(3)] the Kazhdan-Lusztig polynomial associated to $[e, w]$ equals 1;

\item[(4)] $w$ avoids the patterns 4231 and 3412, that is, there do not exits
indices $i_1<i_2<i_3<i_4$ such that the subsequence $w(i_1)w(i_2)w(i_3)w(i_4)$
has the same relative order as 4231 or 3412;
\end{itemize}
see for example Carrell \cite{Car} and Lakshmibai and  Sandhya \cite{Lak}.

When  $X_w$   is smooth,   $P_w(t)$
can be expressed as a product of the factors $1+t+\cdots+t^m$, see Akyildiz and  Carrell \cite{AC}
or Carrell \cite{Car}.
A combinatorial treatment    was given by Gasharov \cite{Gas}.
It is this fact that
motivates us to consider  if the Poincar\'e polynomial of an odd diagram class has
an analogous factorization, as stated in Theorem \ref{main}.

This paper is structured as follows. In Section \ref{Pre}, we collect some
notation, terminology and results used in this paper. In Section \ref{SS3},
we present an algorithm to give a partition of an odd diagram class.
We  prove that the partition is uniform.
Using the results established in Section \ref{SS3},
we finish the proof of Theorem \ref{main} in Section \ref{SS4}.
In Section \ref{SS5}, we discuss problems
concerning  odd diagram classes, including  the self-duality property and
the Kazhdan-Lusztig polynomials of odd diagram classes.

\section{Preliminaries}\label{Pre}

In this section, we  give an overview of the Bruhat order for the symmetric group.
We also describe the legal move operation  introduced by
Brenti,   Carnevale and  Tenner \cite{BCT}, which
plays a  fundamental role in the study of odd diagram classes.

The symmetric group $S_n$ is the Coxeter group of type $A_{n-1}$.
The reflection set is the collection  $\{(i\ j), 1\leq i<j\leq n\}$ of transpositions, and the set
$\{(i\ i+1), 1\leq i\leq n-1\}$ of simple transpositions constitutes  a generating set.
For a permutation  $w\in S_n$, $w\,(i\ j)$ is the permutation obtained by swapping
$w_i$ and $w_j$, while $(i\ j)\,w$ is the permutation obtained by swapping the values
$i$ and $j$. For example, for $w=24513$, we have $w(3\ 4)=24153$ and $(3\ 4)w=23514$.

The Coxeter length $\ell(w)$ of   $w\in S_n$ is the minimum integer $k$
such that $w$ can be expressed as a product of $k$ simple transpositions.
It is well known that $\ell(w)$ equals the number of inversions of $w$:
\begin{equation}\label{length}
\ell(w)=\#\{(w(i), w(j))\colon 1\leq i<j\leq n,\, w(i)>(j)\},
\end{equation}
see for example Bj\"orner and Brenti \cite[Proposition 1.5.2]{BB}.
Notice that $\ell(w)<\ell(w\,(i\ j))$ if and only if $w(i)<w(j)$, and in
this case we denote  $w<w\,(i\ j)$. The transitive closure of
all relations $w<w\,(i\ j)$ forms  the Bruhat order $\leq$ on $S_n$.

Let us recall a combinatorial
 rule for deciding when two permutations are comparable, see   Macdonald \cite[(1.19)]{Mac}.
For two subsets $S, T$ of $[n]$ with the same cardinality,  write
$S\leq T$ if we list the elements of $S$ and $T$ in increasing order,
say $S=\{a_1<a_2<\cdots<a_k\}$ and $T=\{b_1<b_2<\cdots<b_k\}$,
then $a_i\leq b_i$ for each $1\leq i\leq k$.

\begin{prop}\label{Bruhat}
Let $u,v\in S_n$. Then $u\leq  v$ in the Bruhat order if and only if for   $1\leq i\leq n$,
\[\{u(1),u(2),\ldots, u(i)\}\leq \{v(1),v(2),\ldots, v(i)\}.\]
\end{prop}

 For $u,v\in S_n$, we say that
 $u$ is covered by $v$, denoted $u\lhd v$, if there does not exist
 $w\in S_n$ such that $u<w<v$. The covering relation has the following  simple
criterion, see \cite[Lemma 2.1.4]{BB}.

\begin{prop}\label{cover}
Let $u,v\in S_n$. Then $u\lhd v$ if and only if there exist $1\leq i<j\leq n$ such that
$v=u\,(i\ j)$,
$u(i)<u(j)$, and for each $i<k<j$, either $u(k)<u(i)$ or $u(k)>u(j)$.
\end{prop}

Combining \eqref{length} and Proposition \ref{cover},
it follows that   $u\lhd v$ if and only if $v=u\,(i\ j)$ and
 $\ell(v)=\ell(u)+1$.

\begin{figure}[h]
\begin{center}
\begin{tikzpicture}

\def\rectanglepath{-- +(5mm,0mm) -- +(5mm,5mm) -- +(0mm,5mm) -- cycle}

\draw [step=5mm,dotted] (50mm,0mm) grid (70mm,20mm);
\draw (55mm,10mm) \rectanglepath;
\draw (60mm,10mm) \rectanglepath;
\draw (55mm,5mm) \rectanglepath;
\node at (57.5mm,7mm) {*};
\node at (62.5mm,12mm) {*};

\node at (52.5mm,17.5mm) {$\bullet$};
\node at (67.5mm,12.5mm) {$\bullet$};
\node at (62.5mm,7.5mm) {$\bullet$};
\node at (57.5mm,2.5mm) {$\bullet$};

\draw(52.5mm,17.5mm)--(70mm,17.5mm);
\draw(52.5mm,17.5mm)--(52.5mm,0mm);
\draw(67.5mm,12.5mm)--(70mm,12.5mm);
\draw(67.5mm,12.5mm)--(67.5mm,0mm);
\draw(62.5mm,7.5mm)--(70mm,7.5mm);
\draw(62.5mm,7.5mm)--(62.5mm,0mm);
\draw(57.5mm,2.5mm)--(70mm,2.5mm);
\draw(57.5mm,2.5mm)--(57.5mm,0mm);


\node at (2.5mm,17.5mm) {$\bullet$};
\node at (17.5mm,12.5mm) {$\bullet$};
\node at (12.5mm,7.5mm) {$\bullet$};
\node at (7.5mm,2.5mm) {$\bullet$};

\draw(2.5mm,17.5mm)--(20mm,17.5mm);
\draw(2.5mm,17.5mm)--(2.5mm,0mm);
\draw(17.5mm,12.5mm)--(20mm,12.5mm);
\draw(17.5mm,12.5mm)--(17.5mm,0mm);
\draw(12.5mm,7.5mm)--(20mm,7.5mm);
\draw(12.5mm,7.5mm)--(12.5mm,0mm);
\draw(7.5mm,2.5mm)--(20mm,2.5mm);
\draw(7.5mm,2.5mm)--(7.5mm,0mm);

\draw[dotted](50mm,0mm)--(50mm,20mm);
\draw [step=5mm,dotted] (0mm,0mm) grid (20mm,20mm);
\draw (5mm,10mm) \rectanglepath;
\draw (10mm,10mm) \rectanglepath;
\draw (5mm,5mm) \rectanglepath;

\node at (10mm,-5mm) {(a)};\node at (60mm,-5mm) {(b)};

\end{tikzpicture}
\end{center}
\vspace{-6mm}
\caption{(a)   $D(1432)$, \ \ \   (b)  $D_o(1432)$}
\label{RS}
\end{figure}
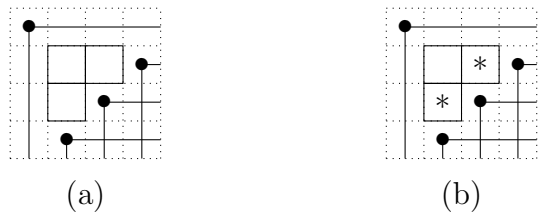

The Rothe diagram $D(w)$ of
$w\in S_n$  is  the   subset
\[D(w)=\{(i,j)\colon w(i)>j,\, i<w^{-1}(j)\} \]
of an $n\times n$ grid.
Alternatively, $D(w)$ can be obtained as follows. For $1\leq i\leq n$, put a dot in the box
$(i, w(i))$, and then delete all boxes lying on the hook with corner at  the box $(i, w(i))$. Then $D(w)$
is exactly the set of the remaining boxes. Figure \ref{RS}(a)
illustrates  the Rothe diagram of $w=1432$.
The odd diagram $D_o(w)$ of $w$ is the subset of $D(w)$ subject to
the parity condition:
\[D_o(w)=\{(i,j)\in D(w)\colon i\not\equiv w^{-1}(j)\ \pmod 2\}.\]
Figure \ref{RS}(b) depicts the odd diagram of $w=1432$, where, as used
in \cite{BCT}, the boxes in $D_o(w)$ are marked with stars.

In the remaining of this section, we give a description of the legal move
operation on odd diagram classes. Assume that $u,v\in S_n$
with $v=u\,(i\ j)$. If $u$ and $v$ have the same odd diagram, then
we say that the transposition $(i\ j)$ is legal for $u$.
The following criterion for a legal transposition  will be used frequently in this paper.

\begin{theo}[\mdseries{Brenti--Carnevale--Tenner \cite[Theorem 4.3]{BCT}}]\label{X2}
Let $u\in S_n$, and $(i\ j)$ be a transposition. Set $m=\min \{u(i), u(j)\}$
and $M=\max \{u(i), u(j)\}$. Then  $(i\ j)$ is a legal transposition
of $u$ if the following conditions are satisfied:
\begin{itemize}
\item[(1)] $i$ and $j$ have the same parity;

\item[(2)] $u(p)<m$ for all $p\in \{i+1, i+3,\ldots, j-1\}$;

\item[(3)] $u(q)\not\in [m, M]$ for all $q\in \{j+1, j+3,\ldots,\}$.

\end{itemize}
\end{theo}

Assume that $u\neq v\in S_n$ have the same odd diagram. It
was observed in \cite{BCT} that
one can apply a legal move to $u$ to obtain a permutation which is
``closer'' to $v$.
Define
\[d(u,v)=\min \{i\colon u^{-1}(i)\neq v^{-1}(i)\}\]
to be the first value  lying at different positions in   $u$ and $v$.
For simplicity, write $k=d(u,v)$. Let $a=u^{-1}(k)$ and
$b=v^{-1}(k)$.

\begin{theo}[\mdseries{Brenti--Carnevale--Tenner \cite[Theorem 4.8]{BCT}}]\label{X1}
The transposition $(a\ b)$ is  legal for $u$.
\end{theo}

Denote $\overline{u}=u\, (a\ b)$. By Theorem \ref{X1}, $\overline{u}$
has the same odd diagram as $u$. Notice
 that $\overline{u}^{-1}(i)=v^{-1}(i)$ for  $i=1,2,\ldots, k$, and so
\[d(\overline{u}, v)>k=d(u,v).\]

Based on Theorem \ref{X1}, it can be shown that in an odd diagram class,
each value   appears in positions with the same parity.

\begin{theo}[\mdseries{Brenti--Carnevale--Tenner \cite[Lemma 6.2]{BCT}}]\label{XX}
Assume that $u, v\in S_n$ have the same odd diagram. Then
\[u^{-1}(i)\equiv v^{-1}(i)\ \pmod 2,\ \ \ \text{for
$1\leq i\leq n$}.\]
\end{theo}

\section{A uniform partition of an odd diagram class}\label{SS3}

Throughout this section, let
\[\mathrm{Perm}(D)=\{w\in S_n\colon D_o(w)=D\}\]
be an odd diagram class in $S_n$ with cardinality greater than 1.
Our aim  is to present a uniform partition
of $\mathrm{Perm}(D)$. In other words, we shall partition
$\mathrm{Perm}(D)$ into blocks with the same cardinality.
Some properties about this uniform partition will be established,
which will be used in the  proof of Theorem
\ref{main} in Section \ref{SS4}.

By Theorem \ref{interval}, $\mathrm{Perm}(D)$ is a  Bruhat
interval, say $[u, v]$. This means that if $w\in S_n$ has odd diagram $D$,
then $u\leq w\leq v$.  By the assumption  $\# \mathrm{Perm}(D)>1$,
we have $u\neq v$.
Fix the following notation
\[k=d(u,v)=\min \{i\colon u^{-1}(i)\neq v^{-1}(i)\}\]
and
\[a=u^{-1}(k),\ \ \ \  b=v^{-1}(k),\]
as used in the preceding section.

\subsection{A partition of $[u, v]$}\label{SB1}

Let us begin with the following lemma.

\begin{lem}\label{L1}
We have $a<b$ and  $u(a)<u(b)$.
\end{lem}

\pf Keep in mind that $u^{-1}(i)=v^{-1}(i)$ for $1\leq i<k$.
Since $u<v$,  Proposition \ref{Bruhat} forces that $a<b$.
Suppose otherwise that $u(a)>u(b)$. By Theorem \ref{X1},
the transposition $(a\ b)$ is legal for $u$. Thus $u\,(a\ b)$ has the same odd diagram
as $u$. The assumption that  $u(a)>u(b)$ leads to $u\,(a\ b)<u$,
contradicting   the minimality of $u$.
\qed

Let
\begin{equation}\label{AS}
\{a=a_1<a_2<\cdots<a_m=b\}=\{a\leq i\leq b\colon u(a)\leq u(i)\leq u(b)\}
\end{equation}
be the set of positions between $a$ and $b$ with values lying in   $[u(a), u(b)]$.
These positions will play a central role in the construction of
the partition of $[u, v]$.

\begin{lem}\label{L2}
The subsequence $u(a_1)u(a_2)\cdots u(a_m)$ of $u$ is   increasing.

\end{lem}

\pf Suppose otherwise there exists $i$ such that $u(a_i)>u(a_{i+1})$.
We claim that the transposition $(a_i\ a_{i+1})$ is legal for $u$.
By Theorem \ref{X1}, the transposition $(a\ b)$ is legal for $u$.
Invoking  Theorem \ref{X2}, we see that
\begin{itemize}
\item[(1)] $a$ and $b$ have the same parity;

\item[(2)] $u(p)<u(a)$ for all $p\in \{a+1, i+3,\ldots, b-1\}$;

\item[(3)] $u(q)\not\in [u(a), u(b)]$ for all $q\in \{b+1, b+3,\ldots\}$.

\end{itemize}

By (2), the positions $a_1, \ldots, a_m$ must have the same parity as $a$,
and so
\begin{equation*}\label{UR1}
a_i\equiv a_{i+1} \pmod 2.
\end{equation*}
It also follows from  (2) that for $p\in \{a_i+1, \ldots, a_{i+1}-1\}$,
\begin{equation*}\label{UR2}
u(p)< u(a_{i+1}).
\end{equation*}
To verify  that $(a_i\ a_{i+1})$ is a legal transposition,
it remains to check that for $q\in \{a_{i+1}+1, a_{i+1}+3,\ldots\}$,
\begin{equation}\label{b}
u(q)\not\in [u(a_{i+1}), u(a_i)].
\end{equation}
This can be seen as follows. Let  $q\in \{a_{i+1}+1, a_{i+1}+3,\ldots\}$. If
$q<b$, then we see from (2) that $u(q)<u(a)$, while  if $q>b$, then
it follows from (3) that $u(q)<u(a)$ or $u(q)>u(b)$. Since $u(a)\leq u(a_i), u(a_{i+1})\leq u(b)$,
we are given  relation \eqref{b}.
This concludes  that the transposition $(a_i\ a_{i+1})$ is  legal for $u$.

By the legality of $(a_i\ a_{i+1})$, the permutation  $u\, (a_i\ a_{i+1})$ has the same odd diagram as $u$.
However, by the assumption  $u(a_i)>u(a_{i+1})$, we are led to $u\, (a_i\ a_{i+1})<u$,
which is contrary to the minimality of $u$. This completes the proof.
\qed

The following lemma  shows  that for any $w\in [u, v]$,  the value $k$
appears in one of the positions $a_1,a_2,\ldots, a_m$.

\begin{lem}\label{L3}
For $w\in [u,v]$, we have
\begin{equation}\label{m}
w^{-1}(k)\in \{a_1,a_2,\ldots, a_m\}.
\end{equation}
\end{lem}

\pf Write $c=w^{-1}(k)$. From Theorem \ref{XX}, it follows that
\begin{equation}\label{UR1}
c\equiv a \equiv b \pmod 2.
\end{equation}
 Since $u\leq w\leq v$, by Proposition \ref{Bruhat},
we are forced that
\begin{equation}\label{p1}
w^{-1}(i)=u^{-1}(i)=v^{-1}(i), \ \ \ \ \text{for  $1\leq i<k$},
\end{equation}
and
\begin{equation}\label{UR11}
a\leq c \leq b.
\end{equation}
By \eqref{UR1} and \eqref{UR11}, we have
\begin{equation}\label{Me}
c\in \{a, a+2,\ldots, b\}.
\end{equation}

In view of the proof of Lemma \ref{L2}, $\{a_1,a_2,\ldots, a_m\}$ is a subset of
$\{a, a+2,\ldots, b\}$.
Suppose to the contrary that
\[c\in \{a, a+2,\ldots, b\}\setminus \{a_1,a_2,\ldots, a_m\}.\]
By the choice of the set $\{a_1,a_2,\ldots, a_m\}$ as given in \eqref{AS}, we  have
either $u(c)<u(a)=k$ or $u(c)>u(b)$.

Case 1. $u(c)<u(a)=k$. In this case, it follows from  \eqref{p1} that $u(c)=w(c)=k$, leading to
a contradiction.

Case 2. $u(c)>u(b)$. According to the values $u(q)$ for $q\in \{b+1, b+3,\ldots,\}$,
the discussion is divided  into   two subcases. Recalling that $(a\ b)$ is a legal transposition for $u$,
it follows from Theorem \ref{X2} that for each $q\in \{b+1, b+3,\ldots,\}$, either $u(q)<u(a)$
or $u(q)>u(b)$.

Subcase 2.1. For each $q\in \{b+1, b+3,\ldots,\}$, either  $u(q)<u(a)$ or $u(q)>u(c)$.
In this case, let us check that $(c\ b)$
is a legal transposition of $u$. By \eqref{Me}, we have $c\equiv b \pmod 2$.
Using again the fact that  $(a\ b)$ is a legal transposition of $u$,
for $p\in \{c+1,\ldots, b-1\}$, we have   $u(p)<u(a)$, and
thus $u(p)<u(b)$. Finally, by the assumption
that $u(q)<u(a)$ or $u(q)>u(c)$ for $q\in \{b+1, b+3,\ldots,\}$, we see that $u(q)\not\in [u(b), u(c)]$.
So the transposition $(c\ b)$ is legal for $u$, and thus  $u\,(c\ b)$ has the same odd diagram
as $u$. However, since $u(c)>u(b)$,
 $u\,(c\ b)$ is smaller than $u$ in the Bruhat order, leading to a contradiction.

Subcase 2.2. There exists $q\in \{b+1, b+3,\ldots,\}$, say $q_0$, such that $u(b)<u(q_0)<u(c)$.
Since $c\not\equiv q_0 \pmod 2$ and $u(c)>u(q_0)$,
the box $(c, u(q_0))$ belongs to   $D_o(u)$. On the other hand, noticing that
\[w(c)=k=u(a)<u(b)<u(q_0),\]
the box $(c, u(q_0))$ cannot belong to $D_o(w)$, contradicting the fact that
$D_o(u)=D_o(w)$. This completes the proof.
\qed

By Lemma \ref{L3}, the interval $[u, v]$ can be partitioned according to
the positions of $k$. Precisely, for $1\leq i\leq m$, set
\[[u, v]^{(i)}=\{w\in [u, v]\colon w^{-1}(k)=a_i\}.\]
To see that each $[u, v]^{(i)}$
is indeed nonempty, we construct a specific permutation $u_i\in S_n$
belonging to $[u, v]^{(i)}$.  As will be seen in Lemma \ref{IO},
$u_i$ is in fact the minimum element of $[u, v]^{(i)}$ in the Bruhat order.

Set $u_1=u$. The constructions of
$u_i$ for $i=2,\ldots, m$ rely on the increasing subsequence $u(a_1)u(a_2)\cdots u(a_m)$.
For $1\leq i\leq m-1$, set
\begin{equation}\label{AAA}
u_{i+1}=u_{i}\,(a_{i}\ a_{i+1}).
\end{equation}

For example, consider the following odd diagram class in $S_9$:
\[[65{\bf4}1{\bf 7}2{\bf 8}3{\bf 9}, 958172634].\]
It is  easily seen that $d(u,v)=4$, $u^{-1}(4)=3$ and $v^{-1}(4)=9$, and so
  $m=4$ and  $(a_1, a_2, a_3, a_4)=(3, 5, 7, 9)$.
The subsequence $u(3)u(5)u(7)u(9)$ is designated  in boldface.
Hence we have
\begin{align*}
&u_1=u=65{\bf4}1{\bf 7}2{\bf 8}3{\bf 9}, \ \ \ \ &u_2=u_1\,(3\ 5)=65{\bf7}1{\bf 4}2{\bf 8}3{\bf 9},\\[5pt]
&u_3=u_2\,(5\ 7)=65{\bf7}1{\bf 8}2{\bf 4}3{\bf 9},\ \ \ \ &u_4=u_3\,(7\ 9)=65{\bf7}1{\bf 8}2{\bf 9}3{\bf 4}.
\end{align*}

\begin{prop}\label{HH}
For $1\leq i\leq m$, the permutation $u_i$ belongs to $[u, v]^{(i)}$.
\end{prop}

\pf It is clear from the construction  that $u_i^{-1}(k)=a_i$.
We still need to show that each $u_i$ has the same odd diagram as $u_1=u$.
To do this, we  assert that for $i=1,\ldots, m-1$, $(a_{i}\ a_{i+1})$
is a legal transposition of $u_{i}$.

Keep in mind that $(a\ b)$ is a legal transposition of $u_1$.
Using  similar arguments as in the proof of Lemma \ref{L2}, we can readily
deduce that
$a_1\equiv a_2 \pmod 2$, $u_1(p)<u(a_1)$ for $p\in \{a_1+1, \ldots, a_2-1\}$, and
$u_1(q)\not\in [u(a_1), u(a_2)]$ for $q\in \{a_2+1, a_2+3, \ldots\}$.
Hence $(a_1\ a_2)$ is a legal transposition of $u_1$.

Analogously, we can verify the assertion for
$i=2,\ldots, m-1$. This implies that  each $u_{i}$ has the
same odd diagram as $u$, and so the proof is complete.
\qed

%

\subsection{The partition is uniform}\label{SB2}

Let us proceed to prove that the partition
\[
[u, v]=\biguplus_{i=1}^m\  [u, v]^{(i)}
\]
is uniform.

\begin{theo}\label{Xmain}
The blocks $[u, v]^{(i)}$ have the same cardinality.
\end{theo}

To give a proof of Theorem \ref{Xmain}, we  construct a bijection
\[\phi_i \colon [u, v]^{(i)} \longrightarrow [u, v]^{(i+1)}\]
 for   $i=1,2,\ldots, m-1$.
 Let $w\in [u, v]^{(i)}$, and set
 \[\phi_i(w)=w\,(a_{i}\ a_{i+1}).\]

Theorem \ref{Xmain} follows from the following proposition.

\begin{prop}\label{PPT}
The map $\phi_i$ is a bijection.
\end{prop}

\pf For $w\in [u, v]^{(i)}$, we first show that  $\phi_i(w)$  belongs to $[u, v]^{(i+1)}$.
It suffices   to verify that  $(a_{i}\ a_{i+1})$
is a legal transposition of $w$.
This can be seen as follows. By Proposition \ref{HH}, the permutation $u_{i+1}$   defined in \eqref{AAA}
belongs to $[u, v]^{(i+1)}$, which,
together with  \eqref{p1}, leads to
\begin{equation*}\label{RR}
d(w, u_{i+1})=k.
\end{equation*}
Applying Theorem \ref{X1} to the pair $w$ and $u_{i+1}$, we see that  $(a_{i}\ a_{i+1})$
is  legal for $w$, and so $\phi_i(w)\in [u, v]^{(i+1)}$.

The reverse construction of $\phi_i$ is clear. Given $w'\in [u, v]^{(i+1)}$,
set $\phi_i^{-1}(w')=w'\,(a_{i}\ a_{i+1})$. Using similar arguments as above,
we can verify   that $\phi_i^{-1}(w')\in [u, v]^{(i)}$. Since $\phi_i^{-1} \circ\phi_i=\phi_i \circ\phi_i^{-1}$
is the identity map,
 $\phi_i$
is a bijection. This completes the proof.
\qed

\begin{figure}[h]
\begin{center}
\begin{tikzpicture}

\node at (40mm,0mm) {$\bullet$};

\node at (40mm,15mm) {$\circ$};
\node at (60mm,15mm) {$\bullet$};
\node at (20mm,15mm) {$\bullet$};

\node at (40mm,30mm) {$\bullet$};
\node at (60mm,30mm) {$\circ$};
\node at (80mm,30mm) {$\bullet$};
\node at (20mm,30mm) {$\circ$};
\node at (0mm,30mm) {$\diamond$};

\node at (40mm,45mm) {$\circ$};
\node at (60mm,45mm) {$\bullet$};
\node at (80mm,45mm) {$\circ$};
\node at (20mm,45mm) {$\diamond$};
\node at (0mm,45mm) {$\diamond$};

\node at (40mm,60mm) {$\diamond$};
\node at (60mm,60mm) {$\circ$};
\node at (20mm,60mm) {$\diamond$};

\node at (40mm,75mm) {$\diamond$};

\draw (40mm,0mm)--(60mm,15mm)--(80mm,30mm)--(80mm,45mm)
--(60mm,60mm)--(40mm,75mm)--(20mm,60mm)--(0mm,45mm)
--(0mm,30mm)--(20mm,15mm)--(40mm,0mm);

\draw (40mm,0mm)--(40mm,15mm)--(60mm,30mm)--(80mm,45mm)
--(40mm,60mm)--(0mm,45mm)--(20mm,30mm)--(40mm,15mm);

\draw (20mm,15mm)--(80mm,30mm)--(60mm,45mm)
--(40mm,60mm)--(20mm,45mm)--(0mm,30mm);

\draw (20mm,15mm)--(20mm,30mm)--(80mm,45mm);

\draw (20mm,15mm)--(40mm,30mm)--(60mm,45mm)--(60mm,60mm);

\draw (40mm,15mm)--(0mm,30mm);

\draw (60mm,15mm)--(60mm,30mm)--(40mm,45mm)
--(20mm,60mm)--(20mm,45mm)--(40mm,30mm)--(60mm,15mm);

\draw (20mm,30mm)--(40mm,45mm)--(60mm,60mm);

\draw (40mm,30mm)--(40mm,45mm);

\draw (60mm,30mm)--(20mm,45mm);
\draw (40mm,60mm)--(40mm,75mm);

\node at (40.5mm,-3mm) {\footnotesize $5431627$};
\node at (40.5mm,12mm) {\footnotesize $5461327$};
\node at (40.5mm,27mm) {\footnotesize $5431627$};
\node at (40.5mm,48mm) {\footnotesize $6471325$};
\node at (40.5mm,63mm) {\footnotesize $7451623$};
\node at (40.5mm,78mm) {\footnotesize $7461523$};

\node at (60.5mm,12mm) {\footnotesize $6431527$};
\node at (60.5mm,27mm) {\footnotesize $6451327$};
\node at (60.5mm,48mm) {\footnotesize $7431625$};
\node at (60.5mm,63mm) {\footnotesize $7461325$};

\node at (20.5mm,12mm) {\footnotesize $5431726$};
\node at (20.5mm,27mm) {\footnotesize $5471326$};
\node at (20.5mm,48mm) {\footnotesize $6451723$};
\node at (20.5mm,63mm) {\footnotesize $6471523$};

\node at (0.5mm,27mm) {\footnotesize $5461723$};
\node at (0.5mm,48mm) {\footnotesize $5471623$};

\node at (80.5mm,27mm) {\footnotesize $7431526$};
\node at (80.5mm,48mm) {\footnotesize $7451326$};

\end{tikzpicture}
\end{center}
\vspace{-6mm}\caption{A uniform partition of an odd diagram class }\label{EE}
\end{figure}
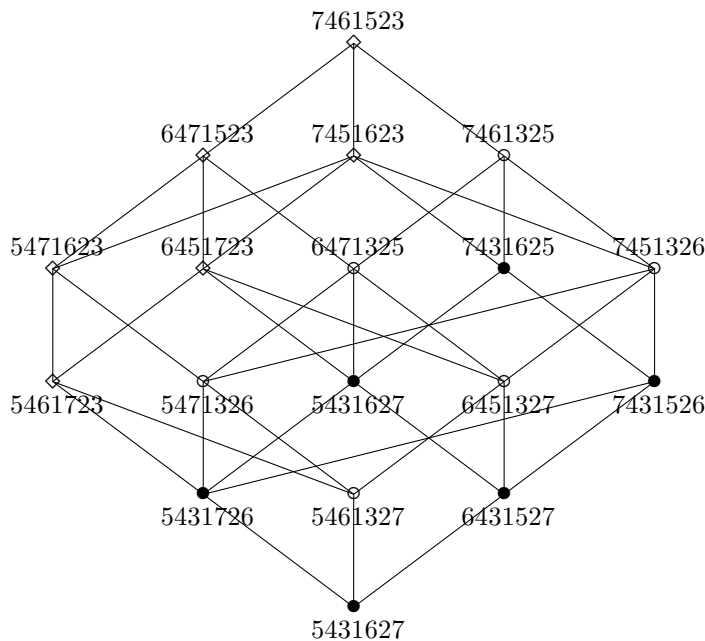

For example, Figure \ref{EE} depicts the following odd diagram class in $S_7$:
\[[u=5431627, v=7461523].\]
We see that $d(u, v)=3$, $u^{-1}(3)=3$ and $v^{-1}(3)=7$. Since
\[\{3\leq i\leq 7\colon u(3)\leq u(i)\leq u(7)\}=\{3, 5,7\},\]
  $[u, v]$ is uniformly  partitioned into the following three blocks:
\begin{align*}
[u, v]^{(1)}&=\{w\in [u, v]\colon w^{-1}(3)=3\},\\
[u, v]^{(2)}&=\{w\in [u, v]\colon w^{-1}(3)=5\},\\
[u, v]^{(3)}&=\{w\in [u, v]\colon w^{-1}(3)=7\},
\end{align*}
which are respectively marked with solid circles, empty circles and diamond symbols.

The bijection $\phi_i$ enjoys the following nice property.

\begin{prop}\label{PP}
For $1\leq i<m$, each permutation  $w\in [u, v]^{(i)}$ is covered by its image $\phi_i(w)$.
\end{prop}

\pf
 For simplicity, write $w'=\phi_i(w)$. Keep in mind that $w^{-1}(k)=a_i$
and $w'^{(-1)}(k)=a_{i+1}$. In view of \eqref{p1}, we see that $d(w, w')=k$, and so we have  $w(a_{i+1})>k$.
By Proposition \ref{cover}, we need to show that
$u(t)\not\in [w(a_i), w(a_{i+1})]$ for $a_i<t<a_{i+1}$.

We use contradiction.
Suppose otherwise that there exists  $a_i<t_0<a_{i+1}$ such that $w(a_i)<w(t_0)<w(a_{i+1})$.
As explained in the proof  of Proposition \ref{PPT}, the transposition  $(a_{i}\ a_{i+1})$
is   legal for $w$. So the conditions in Theorem \ref{X2} are all satisfied
by  $w$ and   $(a_{i}\ a_{i+1})$, from which we can easily check that  $w$ and   $(a_i\ t_0)$
also satisfy the conditions in Theorem \ref{X2}. Thus   $(a_i\ t_0)$
is legal for $w$, implying that   $\overline{w}=w\,(a_i\ t_0)$
belongs to $[u, v]$. However, this would lead to
\[\overline{w}^{-1}(k)=t_0\not\in \{a_1, a_2, \ldots, a_m\},\]
contrary to Lemma \ref{L3}. This completes the proof.
\qed

\subsection{The blocks $[u, v]^{(i)}$ are Bruhat intervals}

In this subsection, we show that each block  in the  partition is a   Bruhat
interval. This will allow us to carry out induction to complete the proof of
Theorem \ref{main} given in Section \ref{SS4}.
Let us first locate  the minimum element of $[u, v]^{(i)}$.

\begin{lem}\label{IO}
For $1\leq i\leq m$, $u_i$ as defined in \eqref{AAA} is the minimum element of $[u, v]^{(i)}$.
\end{lem}

\pf By Proposition \ref{HH}, $u_i$ belongs to $[u, v]^{(i)}$.
It is  clearly true that  $u_1=u$   is the minimum element of $[u, v]^{(1)}$.
We proceed to verify the case for  $u_2$.

Let $w_2\in [u, v]^{(2)}$, and set
\[w_1=\phi_1^{-1}(w_2)\in [u, v]^{(1)}.\]
Then one can find a saturated chain from $u_1$ to $w_1$ in $[u,v]$:
\[u_1=x_1\lhd x_2\lhd \cdots \lhd x_d=w_1. \]
We assert that $x_j\in [u, v]^{(1)}$ for $1\leq j\leq d$.
This can be seen as follows.
By \eqref{p1}, we see that
for $1\leq j\leq d$,
\[x_j^{-1}(s)=u_1^{-1}(s), \ \ \ \ \text{for $1\leq s<k$}.\]
Furthermore, using Proposition \ref{Bruhat} and the fact that
$u_1^{-1}(k)=w_1^{-1}(k)$, we deduce that for $1\leq j\leq d$,
 $x_j^{-1}(k)=u_1^{-1}(k)=a_1$. This verifies  the assertion that $x_j\in [u, v]^{(1)}$.

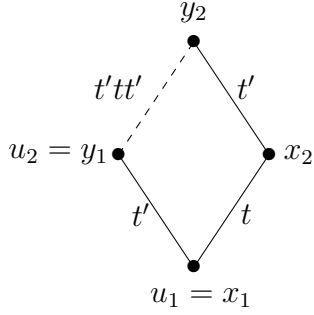
\begin{figure}[h]
\begin{center}
\begin{tikzpicture}

\node at (10mm,30mm) {$\bullet$};
\node at (0mm,15mm) {$\bullet$};
\node at (20mm,15mm) {$\bullet$};
\node at (10mm,0mm) {$\bullet$};
\draw (10mm,0mm)--(20mm,15mm)--(10mm,30mm);
\draw (10mm,0mm)--(0mm,15mm);
\draw [dashed] (0mm,15mm)--(10mm,30mm);

\node at (-8mm,15mm) {$u_2=y_1$};
\node at (10mm,34mm) {$y_2$};
\node at (24mm,15mm) {$x_2$};
\node at (11mm,-4mm) {$u_1=x_1$};

\node at (17mm,7mm) {$t$};\node at (17mm,24mm) {$t'$};
\node at (3mm,7mm) {$t'$};
\node at (0mm,24mm) {$t'tt'$};

\end{tikzpicture}
\end{center}
\vspace{-6mm}
\caption{An illustration of \eqref{OO}}\label{38}
\end{figure}

For $1\leq j\leq d$, define
\[y_j=\phi_1(x_j)\in [u, v]^{(2)}.\]
Note that $y_1=u_2$ and $y_d=w_2$. We claim that $y_1, \ldots, y_d$ form
a saturated  chain:
\[u_2=y_1\lhd y_2\lhd \cdots \lhd y_d=w_2. \]
Let us first check that $y_1\lhd y_2$. Since $x_1\lhd x_2$, we can find a transposition $t$ such
that $x_2=x_1\, t$.
By Proposition \ref{PP}, we have $x_1\lhd y_1$ and
$x_2\lhd y_2$. Write $t'$ for the transposition $(a_1\ a_2)$.
Then   $y_1=x_1\, t'$ and $y_2=x_2\, t'$.
So we have
\begin{equation}\label{OO}
y_2=x_2\, t'=x_1\, t\, t'=x_1\, t'\,(t'\, t\, t')=y_1\,(t'\, t\, t').
\end{equation}
See Figure \ref{38} for an illustration.
Note that $t'\, t\, t'$ is a transposition. Moreover, we have  the
following  length relation
\[\ell(y_2)=\ell(x_2)+1=(\ell(x_1)+1)+1=\ell(y_1)+1,\]
which, along with \eqref{OO}, implies  that $y_1\lhd y_2$.

Using the same arguments as above, we can verify that
$y_j\lhd y_{j+1}$ for $j=2,\ldots, d-1$, and so $y_1, \ldots, y_d$ constitute  a saturated chain.
This yields  that $u_2\leq w_2$, and hence $u_2$ is the minimum of element of $[u, v]^{(2)}$.

To check that $u_3$ is the minimum of element of $[u, v]^{(3)}$, we choose   $w_3\in [u, v]^{(3)}$, and
show that $u_3\leq w_3$ by using completely analogous analysis to the proof for $u_2\leq w_2$.
Continuing the procedure, we eventually  conclude that $u_i$  is
the minimum element of  $[u, v]^{(i)}$ for $1\leq i\leq m$.
This completes the proof.
\qed

We next determine  the maximum element   of $[u, v]^{(i)}$. For $i=m$,
set $v_m=v$. For $1\leq i<m$, set
\[v_i=\phi_i^{-1}(v_{i+1})=v_{i+1}(a_i\ a_{i+1}).\]

\begin{lem}\label{OI}
For $1\leq i\leq m$, $v_i$ is the maximum element of $[u, v]^{(i)}$.
\end{lem}

\pf The proof  is sketched since it is  similar to that of
Lemma \ref{IO}.
Clearly,  $v_m=v$ is the maximum element of $[u, v]^{(m)}$.
We next verify the assertion for  $m-1$.

Let $w\in [u, v]^{(m-1)}$, and set
\[z=\phi_{m-1}(w)\in [u, v]^{(m)}.\]
Locate  a saturated chain from $z$ to $v_m$ in $[u,v]$:
\[z=z_1\lhd z_2\lhd \cdots \lhd z_d=v_m. \]
Similar to   the proof of $x_j\in [u, v]^{(1)}$  in Lemma \ref{IO},
we can verify that   $z_j\in [u, v]^{(m)}$ for $1\leq j\leq d$.

For $1\leq j\leq d$, define
\[w_j=\phi_{m-1}^{-1}(z_j)\in [u, v]^{(m-1)}.\]
Note that $w_1=w$ and $w_d=v_{m-1}$.
We can adopt the same analysis   as used in Lemma \ref{IO}
 to conclude that $w_1\lhd\cdots \lhd w_d$. Hence we have $w\leq v_{m-1}$.
Similar arguments apply to the  case for $i=m-2, \ldots, 1$.
This completes the proof.
\qed

\begin{theo}\label{D1}
For $1\leq i\leq m$, the block $[u, v]^{(i)}$ is the Bruhat interval $[u_i, v_i]$.
\end{theo}

\pf By Lemmas \ref{IO} and \ref{OI}, we have
\[[u, v]^{(i)} \subseteq [u_i, v_i].\]
The reverse inclusion is explained as follows.
For   $w\in [u_i, v_i]$, by \eqref{p1} and the fact $u_i^{-1}(k)=v_i^{-1}(k)=a_i$,
Proposition \ref{Bruhat} forces that
\[w^{-1}(s)=u_i^{-1}(s), \ \ \ \ \text{for $1\leq s\leq k$},\]
which in particular leads to $w^{-1}(k)=a_i$, and so   $w\in [u, v]^{(i)}$,
as desired.
\qed

\section{Proof of Theorem \ref{main}}\label{SS4}

We are now ready to give a proof of Theorem \ref{main}.
Let us start with the observation that the Poincar\'e polynomial of an odd diagram
class $[u,v]$ is the product of the  Poincar\'e polynomial of
the Bruhat interval $[u_m, v]$ and the factor $1+t+\cdots+t^{m-1}$.

\begin{lem}\label{FF}
With the notation in Section \ref{SS3}, we have
\begin{equation}\label{E2}
P_{u, v}(t)=(1+t+\cdots+t^{m-1})P_{u_m, v}(t).
\end{equation}
\end{lem}

\pf By Proposition \ref{PPT} and Theorem \ref{D1},  for $1\leq i\leq m-1$,  $\phi_i$
is a bijection from the interval  $[u_i, v_i]$ to the interval $[u_{i+1}, v_{i+1}]$.
Along with Proposition \ref{PP}  and the fact that $\ell(u_{i+1})=\ell(u_i)+1$,
we obtain that
\begin{align}
 P_{u_i, v_i}(t)&=t^{-\ell(u_i)}\sum_{w\in [u_i, v_i]} t^{\ell(w)}\nonumber\\[5pt]
 &=t^{-\ell(u_i)}t^{-1}\sum_{w\in [u_{i+1}, v_{i+1}]} t^{\ell(w)}\nonumber\\[5pt]
 &=t^{\ell(u_{i+1})-\ell(u_i) -1} P_{u_{i+1}, v_{i+1}}(t)\nonumber\\[5pt]
 &=P_{u_{i+1}, v_{i+1}}(t).\label{E1}
\end{align}
Keep in mind that $u=u_1$ and $v=v_m$.
By \eqref{E1}, we deduce that
\begin{align*}
 P_{u, v}(t)&=t^{-\ell(u_1)}\sum_{i=1}^m \sum_{w\in [u_i, v_i]}t^{\ell(w)}\nonumber\\[5pt]
 &=\sum_{i=1}^m t^{\ell(u_i)-\ell(u_1)} P_{u_i, v_i}(t)\nonumber\\[5pt]
  &=\sum_{i=1}^m t^{i-1}  P_{u_m, v_m}(t)\\[5pt]
  &=(1+t+\cdots+t^{m-1})P_{u_m, v}(t),
\end{align*}
as required.
\qed

We proceed to consider the Bruhat interval $[u_m, v]$.
Carry out the same procedure  in Section \ref{SS3} by replacing $[u, v]$ with $[u_m, v]$. This is sketched as follows.
For convenience, write $w=u_m$.
Set
\[d(w, v)=k',\]
and
\[w^{-1}(k')=a'\ \ \ \text{and }\ \ \ v^{-1}(k')=b'.\]
It should be noted that $k'>d(u, v)=k$ since $w^{-1}(i)=v^{-1}(i)$ for $1\leq i\leq k$.
Let
\[\{a'=a_1'<a_2'<\cdots<a_{m'}'=b'\}=\{a'\leq i\leq b'\colon w(a')\leq w(i)\leq w(b')\}.\]
For $1\leq i\leq m'$, define
\[[w, v]^{(i)}=\{\pi\in [w, v]\colon \pi^{-1}(k')=a_i'\}.\]
Applying the analysis in Subsection  \ref{SB1}, we can show that
$\{[w, v]^{(i)}\colon 1\leq i\leq m'\}$ is a partition of $[w, v]$.

For $1\leq i<m'$, we define a map
\[\phi_i' \colon [w, v]^{(i)} \longrightarrow [w, v]^{(i+1)}\]
by setting
 \[\phi_i'(\pi)=\pi\,(a_{i}'\ a_{i+1}'),  \ \ \text{for $\pi\in [w, v]^{(i)}$}.\]
Similar to Proposition \ref{PPT}, we can show that each $\phi_i'$ is a bijection.
Hence  the blocks
$[w, v]^{(i)}$
 form a uniform partition of
the interval $[w, v]$.

Set $w_1=w$, and
\[w_{i+1}=\phi_i' (w_i), \ \ \ \text{for $i=1,2,\ldots, m'-1$},\]
and set $y_{k'}=v$, and
\[y_{i}=\phi_i'^{-1} (y_{i+1}), \ \ \ \text{for $i=1,2,\ldots, m'-1$}.\]
We can show that the block $[w, v]^{(i)}$ is the Bruhat interval
$[w_i, y_i]$ by using analogous arguments for the proof of Theorem \ref{D1}.

Similar to Lemma \ref{FF}, we deduce that
\[P_{w, v}(t)=(1+t+\cdots+t^{m'-1})P_{w_{m'}, v}(t).\]

Of course, we can continue  the same procedure for the interval $[w_{m'}, v]$.
The procedure eventually terminates, and we reach a proof of Theorem \ref{main}.

\section{Concluding remarks}\label{SS5}

This section is devoted to some observations and problems concerning
odd diagram classes.
As mentioned in Introduction, when a lower interval is rank-symmetric,
its Poincar\'e polynomial factors into polynomials
 of the form $1+t+\cdots+t^m$.
Theorem \ref{main} tells that the Poincar\'e polynomial  of an odd diagram class satisfies such a similar
factorization. It is natural to ask if odd diagram classes  share more    properties satisfied
by rank-symmetric lower intervals.

\subsection{Self-dual odd diagram classes}

The ``top-heavy'' phenomenon  of a lower Bruhat interval $[e, w]$
was  established by Bj\"orner and   Ekedahl \cite{BE}.
 For $0\leq k\leq \ell(w)$,
Let
\[P_k^w=\{u\leq w\colon \ell(u)=k\}\]
denote the rank $k$ component of $[e, w]$.

\begin{theo}[\mdseries{Bj\"orner--Ekedahl \cite{BE}}]\label{HHH}
For $0\leq k\leq \ell(w)/2$,
\begin{equation}\label{TH}
\# P_k^w\leq \# P_{\ell(w)-k}^w.
\end{equation}
\end{theo}

It should be pointed out  that Theorem \ref{HHH} holds in general for parabolic quotients of
Weyl groups.

When the equality in \eqref{TH} holds,  $[e, w]$ is rank-symmetric.
In general, a lower interval is not self-dual.
Gaetz and Gao \cite{GG}  found that the self-duality of $[e, w]$
is determined by  local information of $[e, w]$. Let $\Gamma_w$ (resp., $\Gamma^w$)
denote the bipartite graph on $P_1^w\cup P_2^w$ (resp., $P_{\ell(w)-1}^w\cup P_{\ell(w)-2}^w$)
with edges given by the
covering relations in the Bruhat order.

\begin{theo}[\mdseries{Gaetz--Gao \cite[Theorem 4]{GG}}]\label{UII}
The interval  $[e, w]$ is self-dual if and only if the
bipartite graphs   $\Gamma_w$ and $\Gamma^w$ are isomorphic.
\end{theo}

Note that there are two other criteria for the self-duality of
$[e, w]$ in \cite[Theorem 4]{GG}.

As noticed by Brenti,   Carnevale and  Tenner \cite{BCT},
odd diagram classes are not self-dual in general. For example, the following odd diagram class
\[[654172839, 958172634]\]
is not self-dual. In fact, all odd diagram classes in $S_n$ for $n\leq 8$
are self-dual, and there are 8 and 118 non-self-dual odd diagram classes in $S_9$ and $S_{10}$, respectively.

Obviously, we can define bipartite graphs for   any Bruhat interval.
We checked that the  odd diagrams classes in $S_n$ for $n\leq 10$   satisfy
the bipartite graph criterion for the self-duality,
as given in  Theorem \ref{UII}. Is it possible that
such an observation holds for all odd diagram classes?

\subsection{Kazhdan-Lusztig polynomials}

Let us start with a brief overview of  the Kazhdan-Lusztig polynomials
introduced by Kazhdan and Lusztig \cite{KL}. See
\cite[Chapter 5]{BB} or \cite[Chapter 7]{Hum} for further information.
Let $(W,S)$ be a Coxeter system, and $\leq $ denote the Bruhat order on $W$.
For $w\in W$, a generator $s\in S$ is a  (right) descent if $\ell(ws)<\ell(w)$.
As usual, we use  $D(w)$ to denote the set of descents of $w$.

For $x,y\in W$, the $R$-polynomial  $R_{x,y}(q)$ can be defined in a recursive way:
\begin{itemize}
\item[$\mathrm{(i)}$] $R_{x,y}(q)=0$  if $x\nleq y$;
\item[$\mathrm{(ii)}$] $R_{x,y}(q)=1$ if $x= y$;
\item[$\mathrm{(iii)}$] If $x<y$ and $s\in D(y)$, then
\[R_{x,y}(q)=\left\{
        \begin{array}{ll}
          R_{xs,\,ys}(q), & \hbox{\rm{if} $s\in D(x)$,} \\[5pt]
          qR_{xs,\,ys}(q)+(q-1)R_{x,\,ys}(q), & \hbox{\rm{if} $s \notin  D(x)$.}
        \end{array}
      \right.
\]
\end{itemize}
The Kazhdan-Lusztig polynomials $P_{x, y}(q)$ are the unique family of polynomials
determined by the following
conditions:
\begin{itemize}
\item[$\mathrm{(i)}$] $P_{x,y}(q)=0$  if $x\nleq y$;
\item[$\mathrm{(ii)}$] $P_{x,y}(q)=1$ if $x= y$;
\item[$\mathrm{(iii)}$] if $x\leq y$, then
\[\mathrm{deg}(P_{x, y}(q))\leq \left\lfloor(\ell(y)-\ell(x)-1)/2\right\rfloor\]
and
\[ q^{\ell(y)-\ell(x)}P_{x, y}\left(\frac{1}{q}\right)=\sum_{x\leq z\leq y} R_{x, z}(q) P_{z, y}(q).
\]
\end{itemize}

From the context, no confusion should be caused by the similarity of the notation of the
Kazhdan-Lusztig polynomial $P_{x, y}(q)$ and the Poincar\'e polynomial
$P_{u, v}(t)$.

Recall from Introduction that $[e, w]$ is rank-symmetric if and only if   $P_{e, w}(q)=1$.
We computed that $P_{u, v}(q)=1$ for all odd diagram classes $[u, v]$ in $S_n$ for $n\leq 10$.
It seems reasonable
to conjecture that the Kazhdan-Lusztig polynomial associated to any odd diagram class
is equal to 1.

\noindent
{\bf Remark.} In the lower interval case, $P_{e, w}(q)=1$  if and only if
$[e, w]$ is rank-symmetric.  However, this is not true
for a general interval. For example, let $w_0=n\cdots 2 1$ be the longest permutation in
$S_n$. For any $w\in S_n$, we have $P_{w, w_0}(q)=1$ \cite[(21)]{Bre}, but not every interval  $[w, w_0]$
is rank-symmetric. The following are equivalent conditions  for a general
Kazhdan-Lusztig polynomial equal to 1.

\begin{theo}[\mdseries{Carrell \cite[Theorem C]{Car}}]
Let $(W, S)$ be a Coxeter system, and $T=\cup_{w\in W}wSw^{-1}$ be the
set of reflections. Then the following are equivalent:
\begin{itemize}
\item[(i)] $P_{x, y}(q)=1$;
\item[(ii)] $P_{w, y}(q)=1$ for all $w\in [x, y]$;
\item[(iii)] for all $w\in [x, y]$,
\[\# \{t\in T\colon w<t w\leq y\}=\ell(y)-\ell(w).\]
\end{itemize}
\end{theo}

\vspace{.2cm} \noindent{\bf Acknowledgments.}
This work was
supported by the National Natural
Science Foundation of China  (11971250, 12071320)
and Sichuan Science and Technology Program (Grant No. 2020YJ0006).


\begin{thebibliography}{99}

\bibitem{AC}
E. Akyildiz and J.B. Carrell, A generalization of the Kostant-Macdonald identity,
Proc. Nat. Acad. Sci. U.S.A. 86 (1989),   3934--3937.

\bibitem{BB}
A. Bj\"orner and F. Brenti,
Combinatorics of Coxeter Groups, Grad. Texts in Math., Vol. 231, Springer, New York, 2005.

\bibitem{BE}
 A. Bj\"orner and T.  Ekedahl,   On the shape of Bruhat intervals,
  Ann. Math.  170 (2009), 799--817.

\bibitem{Bre}
F.  Brenti,  Kazhdan--Lusztig and $R$-polynomials from a combinatorial point of view,
 Selected papers in honor of Adriano Garsia (Taormina, 1994),
  Discrete Math. 193 (1998),  93--116.

\bibitem{BC}
F. Brenti and A. Carnevale,
 Proof of a conjecture of Klopsch-Voll on Weyl groups of type $A$, Trans.
Amer. Math. Soc. 369 (2017),  7531--7547.

%


\bibitem{BCT}
F. Brenti, A. Carnevale and B.E. Tenner,
 Odd diagrams, Bruhat order, and pattern avoidance,
 arXiv:2009.08865v1.

\bibitem{Car}
J.B. Carrell,  The Bruhat graph of a Coxeter group, a conjecture of Deodhar, and rational smoothness of Schubert varieties. Algebraic groups and their generalizations: classical methods (University Park, PA, 1991), 53--61, Proc. Sympos. Pure Math., 56, Part 1, Amer. Math. Soc., Providence, RI, 1994.


\bibitem{GG}
C. Gaetz and Y. Gao,  Self-dual intervals in the Bruhat order,
 Selecta Math. (N.S.) 26 (2020),  Paper No. 77, 23 pp.

\bibitem{Gas}
V. Gasharov,  Factoring the Poincar\'e polynomials for the Bruhat order on $S_n$,
 J. Combin. Theory Ser. A 83 (1998),   159--164.

 \bibitem{Hum} J.E. Humphreys, Reflection Groups and Coxeter Groups,  Cambridge
Studies in Advanced Mathematics, No. 29, Cambridge Univ. Press,
Cambridge, 1990.


\bibitem{KL}
D. Kazhdan and G. Lusztig, Representations of Coxeter groups and Hecke algerbas,
Invent. Math. 53 (1979), 165--184.

\bibitem{KV}
B. Klopsch and C. Voll, Igusa-type functions associated to finite formed spaces and their
functional equations, Trans. Amer. Math. Soc. 361 (2009),  4405--4436.


\bibitem{Lak}
V. Lakshmibai and B. Sandhya,   Criterion for smoothness of Schubert varieties in $SL(n)/B$,
 Proc. Indian Acad. Sci. Math. Sci. 100 (1990),  45--52.


\bibitem{Mac}
 I.G. Macdonald, Notes on Schubert Polynomials, Laboratoire de combinatoire et d'informatique math\'ematique (LACIM), Universit\'e du Qu\'ebec \'a Montr\'eal, Montreal, 1991.

\end{thebibliography}
\end{document}